\documentclass{amsart}
\theoremstyle{definition}

\theoremstyle{remark}

\numberwithin{equation}{subsection}

\begin{document}

\title[Jordan counterparts of Rickart and Baer $*$-algebras]
 {Jordan counterparts of Rickart and Baer $*$-algebras}

\author{Shavkat Ayupov}

\address{Institute of Mathematics, National University of Uzbekistan, Tashkent, Uzbekistan}

\email{sh$_-$ayupov@mail.ru}

\author{Farhodjon Arzikulov}

\address{Faculty of Mathematics,
Andizhan State University, Andizhan, Uzbekistan}

\email{arzikulovfn@rambler.ru}

\thanks{This paper partially supported by TWAS, The Abdus Salam, International
Centre, for Theoretical Physics (ICTP), Grant:
13-244RG/MATHS/AS$_-$I-UNESCO FR:3240277696}



\date{}

\dedicatory{}

\commby{}


\begin{abstract}
There are Jordan analogues of annihilators in Jordan algebras which
 are called Jordan annihilators. The present paper is devoted to
investigation of those Jordan algebras every Jordan annihilator of which
is generated by an idempotent as an inner ideal. We prove that
a finite dimensional unital Jordan algebra satisfies this condition
if and only if it has no nilpotent elements, and in this case it is a direct sum of simple Jordan algebras.
\end{abstract}

\maketitle

\section*{Introduction}

It is known that every maximal associative subalgebra of a
finite-dimensional, formally real, unital Jordan algebra over the
field of real numbers is generated by its idempotents and every
family of orthogonal idempotents of this Jordan algebra has the least
upper bound  [1]. Such Jordan algebras are finite
dimensional JB-algebras. In [2] it is proved that the above
conditions are equivalent to the following algebraic condition,
using a Jordan analogies of annihilators: for every subset $S$ of
a JB-algebra $A$, the Jordan annihilator $P^\perp :=\{a\in
A:U_ax=0, (\forall x\in S)\}$ of the set $P=\{a^2: a\in S\}$ is an
inner ideal, generated by an idempotent $e\in A$, i.e.
$P^\perp=U_e(A)$. A Jordan algebra satisfying this algebraic
condition we will call a Baer Jordan algebra.

Present paper is devoted to investigating of Jordan algebras with
more general condition than formally reality. As a corollary of
the Jordan-von Neumann-Wigner theorem [1] we have that
every finite-dimensional, formally real, unital Jordan algebra
over the field of real numbers is a Baer Jordan algebra
[4]. This definition of a Baer Jordan algebra, of course,
is quite compatible to the definition of a Baer $*$-algebra. It
turns out every Baer Jordan algebra has no nilpotent elements. For
a Jordan algebra with no nilpotent elements we also found another
algebraic condition equivalent to the above condition with another
Jordan analogue of the annihilator.

In this paper we also investigate Rickart Jordan algebras which
are Jordan analogies of Rickart $*$-algebras. For the first time,
Rickart Jordan algebras were mentioned in [3]. The
definition of a Rickart Jordan algebra, introduced by us, is also
quite compatible to the definition of a Rickart $*$-algebra.

The last section of the paper is devoted to finite dimensional
Baer Jordan algebras. It is proved that every finite dimensional
unital Jordan algebra without nilpotent elements is a Baer Jordan
algebra and every such Jordan algebra is a direct sum of simple
Jordan algebras.

We use [4], [5] as a standard reference on
notations and terminology.

\begin{center}
{\bf 1. Rickart Jordan algebras}
\end{center}

Let $A$ be a Jordan algebra and
$$
A_+=\{a^2: a\in A\}, \{abc\}=(ab)c+(cb)a-(ac)b,
$$
$$
U_ab=2(ab)a-a^2b, U_a(A)=\{U_ax: x\in A\},
$$
$$
S^\perp :=\{a\in A:U_ax=0, (\forall x\in S)\},
$$
$$
^\perp S:=\{x\in A:U_ax=0, (\forall a\in S)\}
$$
and
$$
^\perp S_+=^\perp S\cap A_+, U_a(A)_+=U_a(A)\cap A_+.
$$
By the Shirshov-Cohn's theorem we have the following identities:
$$
(U_ab)^2=U_aU_ba^2     \,\,\,\,\,\,\,\,      (1)
$$
$$
\{ab(1-a)\}^2=1/4[\{aU_b((1-a)a)(1-a)\}+U_aU_b(1-a)^2+U_{1-a}U_ba^2]
\,\,\,\,\,\,\,\, (2)
$$
By the Macdonald's theorem we have
$$
U_aU_bU_ac=U_{U_ab}c,\,\,\,\,\,\,\,\,\, (3)
$$
$$
4(xy)^2=2xU_yx+U_xy2+U_yx2 \,\,\,\,\,\,\,\,\, (3')
$$
and
$$
U_aU_ab=U_{a^2}b, \,\,\,\,\,\,\,\,\,\,   (4)
$$
if $a$ and $b$ belong to one strongly associative subalgebra of
$A$, then
$$
U_aU_bc=U_{ab}c. \,\,\,\,\,\,\,\,\,\,   (4')
$$
Recall that a subalgebra $B$ of a Jordan algebra $A$ is said to be
strongly associative, if
$$
a(cb)=(ac)b
$$
for all $a\in B$, $b\in B$ and $c\in A$.

Let $A$ be a Jordan algebra. Fix the following conditions:

(A1) for every element $x\in A_+$ there exists an idempotent $e\in
A$ such that $\{x\}^\perp =U_e(A)$;

(B1) for every element $x\in A$ there exists an idempotent $e\in
A$ such that $^\perp \{x\}_+=U_e(A)_+$.

{\it Definition.} A Jordan algebra satisfying condition (A1) is
called {\it a Rickart Jordan algebra}.

For the first time Rickart Jordan algebras have been mentioned in
[3].

{\it Example.} The following Jordan algebra is an example of an
exceptional Jordan algebra, which is a Rickart Jordan algebra: let
$H_3({\Bbb O})$ be the 27-dimensional Jordan algebra of $3\times
3$ hermitian matrices over octonions or Cayley numbers ${\Bbb O}$.
By proposition 3.4 (see section 3) $H_3({\Bbb O})$ is a Rickart
Jordan algebra. So the set $\mathcal{R}$ of all infinite sequences
with components from $H_3({\Bbb O})$ and finite quantity of
nonzero elements is a Rickart Jordan algebra with respect to
componentwise algebraic operations.

Let $A$ be a $*$-algebra, and $S$ is a nonempty subset of $A$, we
write
$$
L(S)=\{x\in A: xs=0 (\forall s\in S)\},
$$
and call $L(S)$ the left-annihilator of $S$. Similarly,
$$
R(S)=\{x\in A: sx=0 (\forall s\in S)\}
$$
denotes the right-annihilator of S. Let $x\in A$. Say
$L(\{x\})=A(1-e)$, $R(\{x\})=(1-f)A$, where $e$ and $f$ are
idempotents.

{\it Definition.} Let $A$ be a $*$-algebra. The involution $*$ is
said to be proper, if $a\neq 0$ implies $aa^*\neq 0$ for any $a\in
A$.

Let $A$ be a $*$-algebra and $A_{sa}=\{a\in A: a^*=a\}$,
$A_+=\{aa^*: a\in A\}$, $U_ab=aba$ and $U_a(A)=\{U_ax: x\in A\}$,
$U_a(A)_+=U_a(A)\cap A_+$. Also we introduce the following
designations:
$$
S^\perp :=\{a\in A:U_ax=0,x\in S\}, ^\perp S:=\{x\in A:U_ax=0,
a\in S\},
$$
$$
^\perp S_+=^\perp S\cap A_+.
$$
Let $A$ be a $*$-algebra. Then the set $A_{sa}$ is a Jordan
algebra with Jordan multiplication $a\cdot b=1/2(ab+ba)$.

{\it Definition.} A Rickart *-algebra is a *-algebra $A$ such
that, for each $x\in A$, $R(\{x\})=gA$ with $g$ a projection
($g^*=g$, $g^2=g$) (note that such a projection is unique). It
follows that $L(\{x\})=(\{x^*\}^r)^*=(hA)^*=Ah$ for a suitable
projection $h$.

{\bf Proposition 1.1.} {\it Let $A$ be a Rickart $*$-algebra. Then
the Jordan algebra $A_{sa}$ satisfies the conditions (A1) and
(B1).}

{\bf Proof.} (A1): Let $x\in A_+$. Then there exists a projection
$e\in A$ such that $L(\{x\})=Ae$ and $axa=0$ for any $a\in
L(\{x\})$. In particular, for all $a\in L(\{x\})_{sa}$ we have
$axa=0$. Since $a=ae=a^*=ea=eae$, $a\in L(\{x\})_{sa}$ we obtain
$$
L(\{x\})_{sa}=eA_{sa}e, eA_{sa}e\subseteq \{x\}^\perp\cap A_{sa}.
\,\,\,\,\,\,\,\,\,\,\,(1)
$$

Now, let $a\in \{x\}^\perp\cap  A_{sa}$. Then $axa=0$. Suppose
$a\notin L(\{x\})$; then $ax\neq 0$, i.e. $ab\neq 0$ for an
element $b\in A_{sa}$ such that $x=bb^*$. Otherwise, if $ab=0$
then $ax=abb^*=(ab)b^*=0b^*=0$. The involution $*$ of a Rickart
$*$-algebra is proper [6]. So $(ab)(ab)^*\neq 0$ and
$(ab)(ab)^*=abb^*a$. Therefore $abb^*a\neq 0$. Hence $axa\neq 0$
and $a\notin \{x\}^\perp\cap A_{sa}$. This is a contradiction.
Hence $a\in L(\{x\})$. Therefore $\{x\}^\perp\cap A_{sa}\subseteq
L(\{x\})_{sa}$. Thus
$$
\{x\}^\perp\cap A_{sa}=L(\{x\})_{sa}=eA_{sa}e=U_e(A_{sa}).
$$

(B1): Take $x\in A_{sa}$. We have $L(\{x\})=Ae$ for some projection
$e$ in $A$, and $(xa)x=xax=0$ for every $a\in
L(\{x\})_+$. Hence
$$
L(\{x\})_+\subseteq ^\perp \{x\}_+.
$$

Now, let $a\in ^\perp \{x\}_+$. Then $xax=0$, i.e. $xbbx=0$, where
$b$ is an element in $A_{sa}$ such that $a=b^2$. We have
$xb(xb)^*=xbbx$. Since the involution $*$ is proper it follows
that $xb=0$. Hence $xa=xbb=0$. Therefore $a\in L(\{x\})_+$. Hence
$^\perp \{x\}_+\subseteq L(\{x\})_+$. Thus
$$
^\perp \{x\}_+=L(\{x\})_+=eA_+e=U_e(A)_+
$$
by (1). $\triangleright$

{\bf Lemma 1.2.} {\it Let $A$ be a $*$-algebra. Suppose the Jordan
algebra $A_{sa}$ of all self-adjoint elements of $A$ is a Rickart
Jordan algebra. Then $A$ has a unit element.}

{\bf Proof.} Take $x=0$. We have $\{x\}^\perp\cap
A_{sa}=U_e(A_{sa})$ for some projection $e$ in $A$. But
$\{x\}^\perp\cap A_{sa}=A_{sa}$. Hence
$U_e(A_{sa})=eA_{sa}e=A_{sa}$. Then $ea=e(ebe)=ebe=a$, where
$a=ebe$, and $ae=(ebe)e=ebe=a$ for every $a\in A_{sa}$. Since
$A=A_{sa}+iA_{sa}$ we have $ea=ae=a$ for all $a\in A$, i.e. $e$ is
an identity element of $A$. $\triangleright$

{\bf Proposition 1.3.} {\it Let $A$ be a $*$-algebra. Suppose the
Jordan algebra $A_{sa}$ of all self-adjoint elements of $A$ is a
Rickart Jordan algebra and the involution $*$ of $A$ is proper.
Then $A$ is a Rickart $*$-algebra.}

{\bf Proof.}  Let $x\in A$. Then $\{xx^*\}^\perp\cap
A_{sa}=U_e(A_{sa})$ for some projection $e$ in $A$. We have
$axx^*a=0$ for all $a\in \{xx^*\}^\perp$. Hence, since the
involution $*$ is proper we have $ax=0$ for all $a\in
\{xx^*\}^\perp\cap A_{sa}$. Therefore
$$
\{xx^*\}^\perp\subseteq A_{sa}\cap L(\{x\})
$$
and
$$
U_e(A)\subseteq  L(\{x\}).
$$
By lemma 1.2 $A$ has a unit element. Hence $e\in L(\{x\})$, i.e.
$ex=0$. Then for each $a\in A$, $aex=0$. Hence $Ae\subseteq
L(\{x\})$.

Let $a\in L(\{x\})$. Then $ax=0$. Hence $a^*axx^*a^*a=0$ and
$a^*a\in \{xx^*\}^\perp\cap A_{sa}$. Therefore $a^*a=ea^*ae$ and
$(1-e)a^*a(1-e)=0$. Hence, since the involution $*$ is proper we
have $a(1-e)=0$, i.e. $a=ae$ and $a\in Ae$. Thus $L(\{x\})=Ae$.
$\triangleright$

{\bf Lemma 1.4.} 1) {\it Let $A$ be a Rickart Jordan algebra. Then

a) $A$ has a unit element;

b) the algebra $A$ has no nilpotent elements. }

{\bf Proof.}  a). Take $x=0$. We have $\{x\}^\perp=U_e(A)$ for
some idempotent $e$ in $A$. But $\{x\}^\perp=A$. Hence $U_e(A)=A$.
Then $ea=eU_eb=U_eb=a$ and $ae=(U_eb)e=U_eb=a$ for all $a\in A$.
Hence we have $ea=ae=a$ for all $a\in A$, i.e. $e$ is a unit
element of $A$.

b). $A$ has a unit element by a) of lemma 1.4. Suppose $a^2=0$,
$a\in A$; then $U_a1=0$ and $a\in \{1\}^\perp$. By the condition
there exists an idempotent $e\in A$ such that $\{1\}^\perp
=U_e(A)$. Since $U_ee=e$ we have $e\in U_e(A)$ and $e\in
\{1\}^\perp$. Hence $U_e1=0$, i.e. $e=0$. Therefore $\{1\}^\perp
=\{0\}$ and $a=0$. Suppose $a^3=0$, $a\in A$; then $a^4=0$ and by
the previous part of the proof $a^2=0$ since $a^4=(a^2)^2=0$.
Hence $a=0$.

Fix $n$. Suppose, if $a^k=0$,  $k\leq n$, $a\in A_+$ then $a=0$
and $a^{n+1}=0$. Then, if $n=2m$ then $a^{2m+2}=(a^{m+1})^2=0$. By
the previous part of the proof $a^{m+1}=0$ and by inductive
supposition $a=0$. if $n=2m+1$ then
$a^{n+1}=a^{2m+2}=(a^{m+1})^2=0$. Again by the previous part of
the proof $a^{m+1}=0$ and $a=0$. Hence by induction, if $a^n=0$
then $a=0$ for every natural number $n$. $\triangleright$

{\bf Lemma 1.5.} {\it Let $A$ be a Rickart Jordan algebra, $x$,
$y\in A_+$. Then the following conditions are equivalent: 1)
$U_xy=0$; 2) $U_yx=0$.}

{\bf Proof.} Suppose $U_xy=0$. Then $U_zU_xz^2=0$ and $(U_zx)^2=0$
by identity (1), where $y=z^2$, $z\in A$. Hence $U_zx=0$ by b) of
lemma 1.4. Therefore
$$
U_yx=U_{z^2}x=U_zU_zx=0
$$
by identity (4).

If $U_yx=0$, then similarly we have $U_xy=0$. $\triangleright$

{\bf Theorem 1.6.} {\it Let $A$ be a Jordan algebra. Then
(A1)$\Longrightarrow$(B1).}

{\bf Proof.} Let $x$ be an element in $A_+$. Then there exists an
idempotent $e\in A$ such that $\{x\}^\perp=U_e(A)$. A has a unit
element by a) of lemma 1.4. Clearly $x\in \{eA(1-e)\}\oplus
U_{1-e}(A)$. Let $x=y^2$, $y\in A$. Then the equality
$$
y=U_ey+\{ey(1-e)\}+U_{1-e}y
$$
gives
$$
y^2=(U_ey)^2+U_ey\{ey(1-e)\}+\{ey(1-e)\}U_ey+\{ey(1-e)\}^2+
$$
$$
\{ey(1-e)\}U_{1-e}y+(U_{1-e}y)^2+U_{1-e}y\{ey(1-e)\}.
\,\,\,\,\,\,\,\,\, (5)
$$
Here,
$$
\{ey(1-e)\}^2=1/4(U_eU_y(1-e)+U_{1-e}U_ye) \,\,\,\,\,\,  (6)
$$
by identity (2). Hence
$$
(U_ey)^2+1/4U_eU_y(1-e)=0
$$
by equalities (5), (6) and $U_ex=0$. But
$$
(U_ey)^2+1/4U_eU_y(1-e)=3/4(U_ey)^2+1/4U_ey^2=3/4(U_ey)^2=0,
$$
and $U_ey=0$ by b) of lemma 1.4. Then $y=U_{1-e}y+\{ey(1-e)\}$. We
will prove that $y\in U_{1-e}(A)$. We have
$$
U_eU_y(1-e)=U_eU_y1-U_eU_ye=U_ey^2-(U_ey)^2=0,
(U_ye)^2=U_yU_ey^2=0
$$
by identity (1), and $U_ye=0$ by b) of lemma 1.4. Hence
$U_{1-e}U_ye=0$ and $\{ey(1-e)\}=0$ by (6) and lemma 1.4.
Therefore $y\in U_{1-e}(A)$ and $x\in U_{1-e}(A)$. Then  $^\perp
\{x\}_+\supset U_e(A)_+$.

Now we will prove that $^\perp \{x\}_+=U_e(A)_+$. Take an element
$c\in ^\perp \{x\}_+$, i.e. $U_xc=0$. Then $U_cx=0$ by lemma 1.5.
Hence $c\in \{x\}^\perp=U_e(A)$. Therefore
$$
^\perp \{x\}_+\subseteq U_e(A)
$$
and
$$
^\perp \{x\}_+=U_e(A)_+.
$$
Let $x$ be an arbitrary element in $A$. Then there exists an
idempotent $e\in A$ such that $^\perp \{x^2\}_+=U_e(A)_+$. By the
previous part of the proof $x\in U_{1-e}(A)$. Hence $^\perp
\{x\}_+\supseteq U_e(A)_+$. Take an element $c\in ^\perp\{x\}_+$,
i.e. $U_xc=0$. Then $U_{x^2}c=U_xU_xc=0$, $U_dU_{x^2}d^2=0$ and
$(U_dx^2)^2=0$ by identities (1) and (4), where $c=d^2$, $d\in A$.
Hence $U_dx^2=0$ by b) of lemma 1.4. Therefore
$U_cx^2=U_{d^2}x^2=U_dU_dx^2=0$ by identity (4). Hence $c\in
\{x^2\}^\perp=U_f(A)$ for some idempotent $f\in A$. But, by the
previous part of the proof $^\perp \{x^2\}_+=U_f(A)_+$. At the
same time $^\perp \{x^2\}_+=U_e(A)_+$. Hence $U_f(A)_+=U_e(A)_+$
and $f=e$. Therefore $c\in U_e(A)_+$ and $^\perp
\{x\}_+\subseteq\mathbb{} U_e(A)+$. Thus
$$
^\perp \{x\}_+=U_e(A)_+.
$$
$\triangleright$

{\bf Theorem 1.7.} {\it Let $A$ be a Jordan algebra with a unit
element. If $A$ has no nilpotent elements, then
(A1)$\Longleftrightarrow$(B1).}

{\bf Proof.} (A1)$\Longrightarrow$(B1) follows by theorem 1.6.

(B1)$\Longrightarrow$(A1): Let $x$ be an element in $A_+$. Then
there exists an idempotent $e\in A$ such that $^\perp
\{x\}_+=U_e(A)_+$. In particular, $U_xe=0$. At the same time
$(U_ex)^2=U_eU_xe=0$ by identity (1). Since $A$ has no nilpotent
elements we obtain $U_ex=0$ and $e\in \{x\}^\perp$. We have
$$
U_{U_ea}x=U_eU_aU_ex=0,
$$
for each $a\in A$, by identity (3). Therefore $U_e(A)\subseteq
\{x\}^\perp$.

Now we will prove that $\{x\}^\perp\subseteq U_e(A)$. Let $a$ be
an element in $\{x\}^\perp$. Then $U_ax=0$. Hence
$U_{a^2}x=U_aU_ax=0$ by identity (4) and
$(U_ya^2)^2=U_yU_{a^2}x=0$ by identity (1), where $y\in A$ such
that $x=y^2$, i.e. $U_ya^2=0$ since $A$ has no nilpotent elements.
Therefore $U_xa^2=U_{y^2}a^2=U_yU_ya^2=0$ by identity (4), i.e.
$a^2\in^\perp \{x\}_+$. Hence $a^2\in U_e(A)_+$. Then similar to
(5) the equality $a=U_ea+\{ea(1-e)\}+U_{1-e}a$ gives
$$
a^2=(U_ea)^2+U_ea\{ea(1-e)\}+\{ea(1-e)\}U_ea+\{ea(1-e)\}^2+
$$
$$
\{ea(1-e)\}U_{1-e}a+(U_{1-e}a)^2+U_{1-e}a\{ea(1-e)\}.
\,\,\,\,\,\,\,\,\, (7)
$$
Hence, since $a^2\in U_e(A)_+$ we have
$$
(U_{1-e}a)^2=0
$$
by equality (7) and by the table of multiplication of Peirce
components. Therefore by the condition of theorem 6
$$
U_{1-e}a=0.
$$
Then $a=U_ea+\{ea(1-e)\}$. We will prove that $a\in U_e(A)$. We
have
$$
\{ea(1-e)\}^2=1/4(U_eU_a(1-e)+U_{1-e}U_ae)
\,\,\,\,\,\,\,\,\,\,\,(8)
$$
by identity (2). But
$$
U_{1-e}U_ae=
U_{1-e}U_a1-U_{1-e}U_a(1-e)=U_{1-e}a^2-(U_{1-e}a)^2=0,
$$
$$
(U_a(1-e))^2=U_aU_{1-e}a^2=0
$$
by identity (1) and $U_a(1-e)=0$ by the condition of theorem 6.
Hence $U_eU_a(1-e)=0$ and $\{ea(1-e)\}=0$. Therefore $a=U_ea$.
Hence $\{x\}^\perp\subseteq U_e(A)$. Thus $\{x\}^\perp=U_e(A)$.
$\triangleright$

{\bf Proposition 1.8.} {\it Let $A$ be a Rickart Jordan algebra,
$x\in A$. There exists a unique idempotent $e$ such that (1)
$xe=x$, and (2) $xy=0$ iff $ey=0$. Explicitly,
$\{x\}^\perp=U_{1-e}A$. The idempotent $e$ is minimal in the
property (1).}

{\bf Proof.} Let $g$ be the idempotent with $\{x\}^\perp=U_gA$,
and set $e=1-g$; clearly $e$ has the property (1). If $h$ is any
idempotent such that $xh=x$, then $(1-h)x=0$, $U_{1-h}x=0$ and
$1-h\in \{x\}^\perp$. Hence $(1-e)(1-h)=1-h$, i.e. $e(1-h)=0$,
$e\leq h$. $\triangleright$

{\it Definition.} Let $A$ be a Rickart Jordan algebra, $x\in A$,
and write $\{x\}^\perp=U_e(A)$ with $e$ idempotent. This
idempotent is unique. Indeed, if $U_e(A)=U_f(A)$ with $e$ and $f$
projections, then $fe=e$ and $ef=f$, so $e=fe=ef=f$. Thus the
idempotent in the definition of a Rickart Jordan algebra is
unique. We will write $e=r(x)$. We call $r(x)$ the range
idempotent of $x$. Thus
$$
\{x\}^\perp=U_{1-r(x)}(A).
$$
For idempotents $e$, $f$ in a Jordan algebra $A$, one writes
$e\leq f$ in case $e\in fAf$ that is, $ef=fe=e$. For projections
$e$, $f$ in a Jordan algebra $A$, the following conditions are
equivalent: $e\leq f$, $e=ef$, $U_e(A)\subseteq U_f(A)$.

{\bf Proposition 1.9.} {\it Let $A$ be a Rickart Jordan algebra
and suppose $\{e_i\}$ is a family of idempotents that has a
supremum $e$. If $x\in A_+$, then $U_ex=0$ iff $U_{e_i}x=0$ for
all $i$.}

{\bf Proof.} In view of proposition 1.8, the following conditions
are equivalent: $U_ex=0$, $r(x)e=0$, $e<1-r(x)$, $e_i<1-r(x)$ for
all $i$, $r(x)e_i=0$ for all $i$, $U_{e_i}x=0$ for all $i$.
$\triangleright$

{\bf Proposition 1.10.} {\it The idempotents of a Rickart Jordan
algebra form a lattice, with
$$
e\vee f = f+r(U1-fe), e\wedge f =e-r(Ue(1-f)).
$$
}

{\bf Proof.} Write $x=U_{1-f}e$ and let $f'=r(x)$. Obviously
$f'\leq f$, so $f+f'$ is an idempotent; we are to show that $f+f'$
serves as $\sup\{e,f\}$. From $x=xf'=(U_{1-f}e)f'=U_{f'}e$ we have
$$
U_{1-f}e=U_{f'}e,
$$
$$
e=U_fe+\{fe(1-f)\}+U_{1-f}e=U_fe+\{fe(1-f)\}+U_{f'}e,
$$
$$
\{fe(1-f)\}^2=U_fU_e(1-f)+U_{1-f}U_ef=
$$
$$
U_fU_e(1-f)-U_{1-f}U_e(1-f)+U_{1-f}U_e1=U_fU_e(1-f)-(U_{1-f}e)^2+U_{1-f}e=
$$
$$
U_fU_e(1-f)-(U_{f'}e)^2+U_{f'}e.
$$
Therefore $(f+f')\{fe(1-f)\}^2=\{fe(1-f)\}^2$. Hence
$$
U_{1-(f+f')}\{fe(1-f)\}^2=0
$$
and
$$
U{\{fe(1-f)\}}U_{1-(f+f')}\{fe(1-f)\}^2=(U_{\{fe(1-f)\}}(1-(f+f')))^2=0,
$$
i.e.
$$
U_{\{fe(1-f)\}}(1-(f+f'))=0
$$
and
$$
U_{1-(f+f')}U_{\{fe(1-f)\}}(1-(f+f'))=(U_{1-(f+f')}\{fe(1-f)\})^2,
$$
$$
U_{1-(f+f')}\{fe(1-f)\}=0.
$$
Therefore $U_{1-(f+f')}e=0$ and $1-(f+f')\in \{e\}^\perp
=U_{1-e}(A)$. Hence $(1-e)(1-(f+f'))=1-(f+f')$ and $e(f+f')=e$,
i.e. $e\leq f+f'$. So $f+f'$ majorizes both $e$ and $f$. Suppose
also $e\leq g$ and $f\leq g$ ($g$ is an idempotent); then
$f=fg=gf$, $e=eg=ge$, so
$$
xg=(U_{1-f}e)g=2(1-f)((1-f)e)g-((1-f)e)g=
$$
$$
2[e-fe-fe-f(fe)]g-(e-fe)g=
$$
$$
2[eg-(fe)g-(fe)g-(f(fe))g]-(eg-(fe)g)=
$$
$$
2[e-f(eg)-f(ge)-f(f(eg))]-(e-f(eg))=
$$
$$
2[e-fe-fe-f(fe)]-(e-fe)=U_{1-f}e=x,
$$
whence $f'\leq g$ and therefore also $f+f'\leq g$. Thus $e\vee f$
exists and is equal to $f+f'$. This establishes the first formula,
and the second follows from it by duality: $e\wedge f$ exists and
$$
e\wedge f=1-[(1-f)\vee (1-e)]=
$$
$$
1-{(1-e)+r(Ue(1-f))}=e-r(Ue(1-f)).
$$
$\triangleright$

{\bf Lemma 1.11.} {\it Let $A$ be a Rickart Jordan algebra, $e$ be
an idempotent. Then $\{e\}^\perp=U_{1-e}(A)$.}

{\bf Proof.} There exists an idempotent $f\in A$ such that
$\{e\}^\perp =U_f(A)$. We have
$$
1-e\in \{e\}^\perp, 1-e\leq f
$$
and
$$
U_fe=0.
$$
Then
$$
\{ef(1-e)\}^2=1/4[\{eU_f((1-e)e)(1-e)\}+U_eU_f(1-e)^2+U_{1-e}U_fe^2]=
$$
$$
1/4[U_eU_f(1-e)^2+U_{1-e}U_fe^2]=
1/4[U_eU_fe+U_ef+U_{1-e}U_fe]=U_ef=0
$$
by formula (2) and lemma 1.5. Hence $\{ef(1-e)\}=0$ and
$fe=f(U_fe+\{ef(1-f)\}+U_{1-f}e)=fU_{1-f}e=0$. Therefore $e\leq
1-f$. Hence $e=1-f$, $1-e=f$ and $\{e\}^\perp=U_{1-e}(A)$.
$\triangleright$

\medskip

\begin{center}
{\bf 2. Baer Jordan algebras}
\end{center}

\medskip

Let $A$ be a Jordan algebra. Fix the following conditions:

(A2) for every subset $S\subseteq A_+$ there exists an idempotent
$e\in A$ such that $S^\perp =U_e(A)$;

(B2) for every subset $S\subseteq A$ there exists an idempotent
$e\in A$ such that $^\perp S_+=U_e(A)_+$.

These conditions are Jordan analogies of Baer conditions for Baer
$*$-algebras.

{\it Definition.} A Jordan algebra satisfying condition (A2) is
called {\it a Baer Jordan algebra}.

{\it Example.} The following Jordan algebra is an example of an
exceptional Jordan algebra, which is a Baer Jordan algebra: let
$H_3({\Bbb O})$ be the 27-dimensional Jordan algebra of $3\times
3$ hermitian matrices over octonions or Cayley numbers ${\Bbb O}$.
By corollary 3.6 (see section 3) $H_3({\Bbb O})$ is a Baer Jordan
algebra. Hence the set $\mathcal{B}$ of all infinite sequences
with components from $H_3({\Bbb O})$ is a Baer Jordan algebra with
respect to componentwise algebraic operations.

Note that the Rickart Jordan algebra $\mathcal{R}$ of all infinite
sequences with components from $H_3({\Bbb O})$ and finite quantity
of nonzero elements is not a Baer Jordan algebra.

{\it Definition.} A Baer $*$-algebra is a $*$-algebra $A$ such
that, for every nonempty subset $S$ of $A$, $R(S)=gA$ for a
suitable projection $g$. (It follows that
$L(S)=(R(S^*))^*=(hA)^*=Ah$ for a suitable projection $h$.)

The following proposition holds:

{\bf Proposition 2.1.} {\it Let $A$ be a Baer $*$-algebra. Then
the Jordan algebra $A_{sa}$ satisfies conditions (A2) and (B2).}

{\bf Proof.} (A2): Let $S\subseteq A_+$. Then there exists a projection
$e\in A$ such that $L(S)=Ae$, where $L(S)=\{a\in A:
as=0, s\in S\}$. Then $asa=0$ for any $a\in L(S)$. In particular
for all $s\in S$ and $a\in L(S)\cap A_{sa}$ we have $asa=0$. Since
$a=ae=a^*=ea=eae$, $a\in L(S)\cap A_{sa}$ it follows that
$$
L(S)\cap A_{sa}=eA_{sa}e, eA_{sa}e\subseteq S^\perp\cap A_{sa}.
\,\,\,\,\,\,\,(9)
$$

Now, let $a\in S^\perp\cap A_{sa}$. Then $asa=0$ for all $s\in S$.
Suppose $a\notin L(S)$; then there exists $s\in S$ such that
$as\neq 0$, i.e. $ab\neq 0$ for an element $b\in A_{sa}$ such that
$s=bb^*$. Otherwise, if $ab=0$ then $as=abb^*=(ab)b^*=0b^*=0$. The
involution $*$ of any Baer $*$-algebra is proper [6]. So
the involution $*$ of the $*$-algebra $A$ is proper and
$(ab)(ab)^*=abb^*a$. Therefore $abb^*a\neq 0$. Hence $asa\neq 0$
and $a\notin S^\perp\cap A_{sa}$. This is a contradiction. Hence
$a\in L(S)$. Therefore $S^\perp\cap A_{sa}\subseteq L(S)_{sa}$.
Thus
$$
S^\perp\cap A_{sa}=L(S)_{sa}=eA_{sa}e=U_e(A_{sa}).
$$

(B2): Take $S\subseteq A_{sa}$. We have $L(S)=Ae$ for some
projection $e$ in $A$. At the same time $sas=s(as)=0$ for all
$s\in S$ and $a\in L(S)\cap A_+$. Hence
$$
L(S)\cap A_+\subseteq  {^\perp S_+}.
$$

Now, let $a\in{^\perp S_+}$. Then $sas=0$ for any $s\in S$. At the
same time $sbbs=0$ for any $s\in S$, where $b$ is an element in
$A_{sa}$ such that $a=b^2$. We have $sb(sb)^*=sbbs$. Since the
involution $*$ is proper it follows that $sb=0$ for all $s\in S$.
Hence $sa=sbb=0$ for every $s\in S$. Therefore $a\in L(S)\cap
A_+$. Hence $^\perp S_+\subseteq L(S)\cap A_+$. Thus
$$
^\perp S_+=L(S)\cap A_+=eA_+e=U_e(A)_+
$$
by (9). $\triangleright$

{\bf Lemma 2.2.} {\it Let $A$ be a $*$-algebra and $A_{sa}$ the
Jordan algebra of all self-adjoint elements of $A$ satisfying the
following condition: for each subset $S\subseteq A_+$ there exists
a projection $e\in A$ such that $S^\perp\cap A_{sa}=U_e(A_{sa})$.
Then $A$ has a unit element.}

{\bf Proof.} This lemma follows by lemma 1.2. $\triangleright$

{\bf Proposition 2.3.} {\it Let $A$ be a $*$-algebra. Suppose the
Jordan algebra $A_{sa}$ of all self-adjoint elements of $A$ is a
Baer Jordan algebra and the involution $*$ of $A$ is proper. Then
$A$ is a Baer $*$-algebra.}

{\bf Proof.}  Let $S\subseteq A$. Then $T^\perp\cap
A_{sa}=U_e(A_{sa})$ for some projection $e$ in $A$, where
$T=\{ss^*:s\in S\}$. We have $ass^*a=0$ for all $a\in T^\perp\cap
A_{sa}$ and $s\in S$. Hence, since the involution $*$ is proper we
have $as=0$ for all $a\in T^\perp\cap A_{sa}$ and $s\in S$.
Therefore
$$
T^\perp\cap A_{sa}\subseteq L(S)
$$
and
$$
U_e(A)\subseteq L(S).
$$
By lemma 2.2 $A$ has a unit element. Hence $e\in L(S)$, i.e.
$es=0$ for any $s\in S$. Then $aes=0$ for every $a\in A$ and $s\in
S$. Hence $Ae\subseteq L(S)$.

Let $a$ be an element of $L(S)$. Then $as=0$ for any $s\in S$.
Hence $a^*ass^*a^*a=0$ and $a^*a\in T^\perp\cap A_{sa}$. Therefore
$a^*a=ea^*ae$ and $(1-e)a^*a(1-e)=0$. Hence, since the involution
$*$ is proper we have $a(1-e)=0$, i.e. $a=ae$ and $a\in Ae$. Thus
$$
L(S)=Ae.
$$
$\triangleright$

{\bf Lemma 2.4.} {\it Let $A$ be a Baer Jordan algebra. Then

a) $A$ has a unit element;

b) the algebra $A$ has no nilpotent elements. }

{\bf Proof.} This lemma follows by lemma 1.4. $\triangleright$

{\bf Theorem 2.5.} {\it Let $A$ be a Jordan algebra. Then
(A2)$\Longrightarrow$(B2).}

{\bf Proof.}  By a) of lemma 2.4 $A$ has a unit element. Let $S$
be a subset of $A$. Then there exists an idempotent $e\in A$ such
that $(^\perp S_+)^\perp=U_e(A)$. But we have
$$
U_{1-e}(A_+)^\perp=U_e(A)
$$
by lemma 1.11 and, if $a\notin U_{1-e}(A)_+$ for some $a\in ^\perp
S_+$, then $U_ea\neq 0$ and the set $(^\perp S_+)^\perp$ does not
equal to $U_e(A)$. Indeed, since $(^\perp S_+)^\perp=U_e(A)$ we
have $U_ea=0$. Let $b\in A$ be an element such that $a=b^2$. The
equality $b=U_eb+\{eb(1-e)\}+U_{1-e}b$ gives
$$
b^2=(U_eb)^2+U_eb\{eb(1-e)\}+\{eb(1-e)\}U_eb+\{eb(1-e)\}^2+
$$
$$
\{eb(1-e)\}U_{1-e}b+(U_{1-e}b)^2+U_{1-e}b\{eb(1-e)\}.
\,\,\,\,\,\,\,\,\, (10)
$$
We have
$$
\{eb(1-e)\}^2=1/4(U_eU_b(1-e)+U_{1-e}U_be)
\,\,\,\,\,\,\,\,\,\,\,\,(11)
$$
by identity (2). Hence, since $U_ea=0$ we have
$$
U_ebU_eb+1/4U_eU_b(1-e)=0
$$
by (10) and by the table of multiplication of Peirce components.
But
$$
U_ebU_eb+1/4U_eU_b(1-e)=3/4U_ebU_eb+1/4 U_eb^2=3/4(U_eb)^2=0.
$$
Hence by b) of lemma 2.4 $U_eb=0$. Then $b=U_{1-e}b+\{eb(1-e)\}$.
We prove that $\{eb(1-e)\}=0$. We have
$$
U_eU_b(1-e)=U_eU_b1-U_eU_be=U_eb^2-(U_eb)^2=0,
$$
$$
(U_be)^2=U_bU_eb^2=0
$$
by identity (1) and $U_be=0$ by b) of lemma 2.4. Hence
$U_{1-e}U_be=0$, $\{eb(1-e)\}^2=0$ by (11) and $\{eb(1-e)\}=0$.
Therefore $b\in U_{1-e}(A)$ and $a\in U_{1-e}(A)$. Thus $^\perp
S_+\subseteq U_{1-e}(A)_+$. By lemma 2.4 $A$ has a unit element.
Therefore $1-e\in ^\perp S_+$, i.e. $U_a(1-e)=0$, for all $a\in
S$. By condition (A2) and by lemma 1.11 $\{1-e\}^\perp=U_e(A)$.

Now we prove that $U_{1-e}(A)_+\subseteq ^\perp S_+$. Since
$(^\perp S_+)^\perp=U_e(A)$ we have $S\subseteq U_e(A)$. Hence
$$
U_sU_{1-e}a=U_{se}U_{1-e}a=U_sU_eU_{1-e}a=U_sU_{e(1-e)}a=U_sU_0a=0
$$
for all $s\in S$ and $a\in A$ by identity ($4'$). Therefore
$U_{1-e}(A)_+\subseteq ^\perp S_+$. Thus $^\perp
S_+=U_{1-e}(A)_+$. $\triangleright$

{\bf Theorem 2.6.} {\it Let $A$ be a Jordan algebra. If $A$ has no
nilpotent elements, then (A2)$\Longleftrightarrow$(B2).}

{\bf Proof.} (A2)$\Longrightarrow$(B2) follows by theorem 2.5.

(B2)$\Longrightarrow$(A2): By 2) of lemma 2.4 $A$ has a unit
element. Let $S$ be a subset of $A_+$. Then there exists an
idempotent $e\in A$ such that $^\perp (S^\perp)_+=U_e(A)_+$. If
$a\notin U_{1-e}(A)$ for some $a\in S^\perp$, then $U_ae\neq 0$
and the set $^\perp(S^\perp)_+$ does not equal to $U_e(A)_+$.
Indeed, since $^\perp(S^\perp)_+=U_e(A)_+$ we have $U_ae=0$. Then
$U_aU_ae=U_{a^2}e=0$ and $U_eU_{a^2}e=0$ and $(U_ea^2)^2=0$ by
identity (4). Hence $U_ea^2=0$ by the condition of theorem 2.6. As
in the proof theorem 2.5 we have $U_ea=0$. Then
$$
a=U_{1-e}a+\{ea(1-e)\}.
$$
Since $a\notin U_{1-e}(A)_+$ we have $\{ea(1-e)\}\neq 0$. Then
$U_ea^2\neq 0$ by the proof of theorem 2.5. This is a
contradiction. Therefore $a\in U_{1-e}(A)$ and $S^\perp\subseteq
U_{1-e}(A)$.

Now we have to show that $U_{1-e}(A)\subseteq S^\perp$. Since
$^\perp (S^\perp)_+=U_e(A)_+$, $S\subseteq ^\perp(S^\perp)_+$ we
obtain
$$
S\subseteq U_e(A)_+
$$
and
$$
U_{1-e}s=0
$$
for all $s\in S$. Therefore $1-e\in S^\perp$. Then for any $a\in
U_{1-e}(A)$ we have $U_as=U_{U_{1-e}a}s=U_{1-e}U_aU_{1-e}s=0$ for
all $s\in S$ by identity (3). Hence $U_{1-e}(A)\subseteq S^\perp$
and $U_{1-e}(A)=S^\perp$. This concludes the proof.
$\triangleright$

{\bf Theorem 2.7.} {\it The following conditions are equivalent:
(a) $A$ is a Baer Jordan algebra; (b) $A$ is a Rickart Jordan
algebra and the set of all idempotents of $A$ is a complete
lattice.}

{\it Proof.}  (a)$\Longrightarrow$(b): By proposition 1.10 of
section 1 all idempotents of a Baer Jordan algebra form a lattice.
Let $\{e_i\}$ be an arbitrary family of idempotents in A. Then
there exists an idempotent $e$ in $A$ such that
$\{e_i\}^\perp=U_{1-e}(A)$. We prove that $\sup_i e_i=e$. Indeed,
let $g$ be an idempotent such that $g\geq e_i$ for all $i$. Then
$e_i(1-g)=0$ for all $i$. Hence $U_{1-g}e_i=0$ for all $i$ and
$1-g\in \{e_i\}^\perp$. Hence $1-g\in U_{1-e}(A)$ and
$(1-e)(1-g)=1-g$, $1-g-e+eg=1-g$, $eg=e$. Therefore $e\leq g$.
Since $g$ is chosen arbitrarily $\sup_i e_i=e$. At the same time
$\inf_i e_i=1-\sup_i (1-e_i)$.

(b)$\Longrightarrow$(a):  Let $S=\{x_\alpha :\alpha\in \Omega\}$
be any subset of $A_+$; we are to show that there exists an
idempotent $e\in A$ such that $S^\perp =U_e(A)$. Now,
$$
S^\perp=\cap_{\alpha\in\Omega} \{x_\alpha\}^\perp.
$$
Since $A$ is a Rickart Jordan algebra, there exists an idempotent
$e_\alpha\in A$ such that $\{x_\alpha\}^\perp=U_{e_\alpha}(A)$ for
any $\alpha$. Let $e=inf_\alpha e_\alpha$ in $A$, which exists by
the hypothesis of the theorem. Then
$$
e=1-sup_\alpha(1-e_\alpha)
$$
in $A$. We have
$$
U_e(A)\subseteq \{x_\alpha\}^\perp
$$
for any $\alpha$. Hence $U_e(A)\subseteq S^\perp$. We prove that
$\cap_{\alpha\in\Omega}\{x_\alpha\}^\perp\subseteq U_e(A)$. Let
$a$ be a positive element in
$\cap_{\alpha\in\Omega}\{x_\alpha\}^\perp$. Then $a\in
U_{e_\alpha}(A)$ and $ae_\alpha=a$, i.e. $U_{1-e_\alpha}a=0$ for
any $\alpha$. Hence
$$
U_{\sup_\alpha (1-e_\alpha)}a=0,
$$
i.e.
$$
U_{1-e}a=0
$$
by proposition 1.9 of section 1. Therefore $U_a(1-e)=0$ by lemma
1.5 of section 1. Hence $a\in U_e(A)$, i.e.
$$
(S^\perp)_+\subseteq U_e(A).
$$

Now, let $a$ be an arbitrary element in
$\cap_{\alpha\in\Omega}\{x_\alpha\}^\perp$. Then $a^2\in U_e(A)$,
i.e. $(1-e)a^2=0$. Hence
$$
U_{1-e}a^2=0, U_aU_{1-e}a^2=(U_a(1-e))^2=0
$$
and $U_a(1-e)=0$. By lemma 1.11 of section 1
$\{1-e\}^\perp=U_e(A)$ and $a\in \{1-e\}^\perp$. Therefore $a\in
U_e(A)$. Thus $S^\perp=U_e(A)$ and $A$ is a Baer Jordan algebra.
$\triangleright$

\medskip

\begin{center}
{\bf 3. Finite - dimensional Baer Jordan algebras}
\end{center}

\medskip

By lemma 2.9.4 in [4] and by its proof we have the
following lemma:

{\bf Lemma 3.1.} {\it Let $A$ be a finite - dimensional, unital
Jordan algebra over $\Bbb R$ with no nilpotent elements. Then we
have:

(i) An idempotent $p$ in $A$ is minimal if and only if
$\{pAp\}=Fp$, , where $F={\Bbb R}$ or ${\Bbb C}$.

(ii) Any element of $A$ is contained in some maximal associative
subalgebra of $A$, and every such subalgebra is of the form
$F_1p_1\oplus F_2p_2\oplus\dots \oplus F_np_n$, where
$p_1$,$\dots,$ $p_n$ are pairwise orthogonal minimal idempotents
with sum $1$, $F_i\in \{{\Bbb R}, {\Bbb C}\}$, $i=1,2,\dots,n$.

(iii) If $p$ and $q$ are orthogonal minimal idempotents in $A$ and
$a\in \{pAq\}$, then $a^2=\lambda (p+q)$, where $\lambda\in {\Bbb
R}$ or ${\Bbb C}$.}

{\it Definition.} Let $A$ be a Jordan algebra and $p$, $q$ be
minimal idempotents. We say that $p$, $q$ are {\it connected}, if
$\{pAq\}\neq \{0\}$.

{\bf Lemma 3.2.} {\it Let $A$ be a finite-dimensional, unital
Jordan algebra. Let $p_1, p_2, \dots,p_n$ be a family of pairwise
orthogonal minimal idempotents with sum 1. Write $p_i\sim p_j$ if
$p_i$ and $p_j$ are connected. Then $\sim$ is an equivalence
relation, and for all $i$, $\sum_{p_j\sim p_i}p_j$ is a central
idempotent in $A$.}

{\bf Proof.} It is clear that the relation $\sim$ is reflexive and
symmetric. We prove transitivity of this relation.

If $p_i\sim p_j\sim p_k$ then there exist nonzero elements $a$ and
$b$ in the Peirce components $\{p_iAp_j\}$, $\{p_jAp_k\}$
respectively. Let $\lambda$, $\mu$ be numbers such that
$a^2=\lambda(p_i+p_j)$, $b^2=\mu(p_j+p_k)$ respectively. We have
$ab\in \{p_iAp_j\}$ and
$$
4(ab)^2=2aU_ba+U_ab^2+U_ba^2
$$
by identity ($3'$). Hence
$$
U_ab^2=U_a(\mu(p_j+p_k))=\mu a^2-\mu
a^2(p_j+p_k)=\lambda\mu(p_i+p_j)-\lambda\mu p_j,
$$
$$
U_ba^2=U_b(\lambda(p_i+p_j))=\lambda b^2-\lambda
b^2(p_j+p_i)=\lambda\mu(p_j+p_k)-\lambda\mu p_j,
$$
$$
2aU_ba=2a(2(b(ba)-b^2a)=-2a\mu((p_j+p_k)a)=-\mu
a^2=-\lambda\mu(p_i+p_j)
$$
and
$$
4(ab)^2=2aU_ba+U_ab^2+U_ba^2=
$$
$$
-\lambda\mu(p_i+p_j)+\lambda\mu(p_i+p_j)-\lambda\mu
p_j+\lambda\mu(p_j+p_k)-\lambda\mu p_j=
$$
$$
\lambda\mu(p_j+p_k)-2\lambda\mu p_j=\lambda\mu(p_k-p_j)\neq 0.
$$
So $ab\neq 0$. Therefore $p_i\sim p_k$
and the given relation is transitive, so $\sim$ is an equivalence
relation.

Let $e=\sum_{p_j\sim p_i}p_j$. To show that $e$ is central it is
enough to show that $e$ operator commutes with any element of
$A_{kl}=\{p_kAp_l\}$ for any $k$, $l$. This is clear if $p_k\sim
p_i\sim p_l$, or if neither $p_k\sim p_i$ nor $p_l\sim p_i$, for
then, if $a\in A_{kl}$, $a\circ e=a$ or $a\circ e=0$, so 2.5.5 in
[4] shows $[T_a,T_e]=0$. If, on the other hand, $p_k\sim
p_i$ and the projections $p_i$, $p_l$ are not equivalent, then
$A_{kl}=0$ since $\sim$ is transitive. This completes the proof.
$\triangleright$

{\bf Theorem 3.3.} {\it Every finite dimensional unital Jordan
algebra $A$ without nilpotent elements is a direct sum of simple
algebras. If $A$ is simple then it contains $n\geq 1$ pairwise
orthogonal and connected minimal idempotents with sum 1.}

{\bf Proof.} This theorem can be proven as in the proof of theorem
2.9.6 in [4].

Indeed, by lemma 3.1 $A$ contains $n\geq 1$ pairwise orthogonal
minimal idempotents $p_1$, $\dots$, $p_n$ with sum 1. If $e$ is a
central projection in $A$ then $eA$ is an ideal, and $A=eA+(1-e)A$
is an algebra direct sum. Therefore, from lemma 3.2 we may
conclude two things. First, if $A$ is simple there can be no
nontrivial central projections, and so $p_1$, $\dots$, $p_n$ are
all pairwise connected. Secondly, $A$ is a direct sum of algebras
which do contain pairwise orthogonal and connected minimal
idempotents with sum 1. If $A$ satisfies this requirement we must
show that $A$ is simple. So let $I$ be an ideal in $A$, let $a\in
I$ be nonzero, and write $a=\sum_{i\leq j}a_{ij}$, where
$a_{ij}\in A_{ij}=\{p_iAp_j\}$. Pick $i$, $j$ so that $a_{ij}\neq
0$. If $i=j$, then $a_{ii}=\{p_iap_i\}\in I$. Since by (i) of
lemma 3.1 $a_{ii}=\lambda p_i$ with $\lambda\neq 0$, $p_i\in I$.
If $i\neq j$ then by (iii) of lemma 3.1
$a_{ij}^2=\lambda(p_i+p_j)$, $a_{ij}=2\{p_iap_j\}\in I$, so again
$p_i=\lambda^{-1}p_ia_{ij}^2$ belongs to $I$. This implies that
any $b\in A_{ki}$ belongs to $I$, for any $k$, since $b=2p_ib$ (or
$b=p_ib$, if $k=i$). Repeat the previous argument with any nonzero
element of $A_{ki}$ to conclude $p_k\in I$, so $1=\sum p_k \in I$,
and thus $I=A$, i.e. $A$ is simple. $\triangleright$

{\bf Proposition 3.4.} {\it Every finite dimensional unital Jordan
algebra $A$ without nilpotent elements is a Rickart Jordan
algebra.}

{\bf Proof.} Let $a\in A$ be an arbitrary positive element. There
exists a maximal associative subalgebra $\mathcal{A}_o$ containing
$a$ and
$$
\mathcal{A}_o=F_1p_1\oplus F_2p_2\oplus\dots F_np_n,
$$
where $p_1$, $\dots$,$p_n$ are pairwise orthogonal minimal
idempotents with sum 1, $F_i\in \{{\Bbb R}, {\Bbb C}\}$,
$i=1,2,\dots,n$. Hence there exist $\lambda_i\in F_i$,
$i=1,2,\dots,k$, $k\leq n$, such that
$$
a=\lambda_1e_1+\lambda_2e_2+\dots+\lambda_ke_k
$$
for some subset $\{e_1,e_2,\dots,e_k\}\subseteq
\{p_1,p_2,\dots,p_n\}$.

We prove that $\{a\}^\perp=U_{1-e}(A)$, where
$e=e_1+e_2+\dots+e_k$.

Indeed, it is clear that $U_{1-e}(A)\subseteq \{a\}^\perp$. Let
$b$ be an element in $\{a\}^\perp$. Then $U_ba=0$. There exist
elements $c$, $d$ in $A$ such that $cb$, $da$ are idempotents,
$e=da$ and
$$
U_cU_ba=U_{cb}a=0.
$$
Hence $U_a(cb)=0$ by lemma 1.5 of section 1 and
$$
U_dU_a(cb)=U_{da}(cb)=U_e(cb)=0
$$
by ($4'$). Let $f=cb$. Then $U_ef=0$. Therefore $U_fe=0$. Hence
$$
4(ef)^2=2eU_fe+U_ef^2+U_fe^2=
$$
$$
2eU_fe+U_ef+U_fe=0.
$$
by Macdonald's identity ($3'$). Hence $fe=0$. Thus $f\leq 1-e$ and
$b\in U_{1-e}(A)$. This completes the proof. $\triangleright$

{\bf Proposition 3.5.} {\it The set of all idempotents of every
finite dimensional unital Jordan algebra $A$ without nilpotent
elements is a complete lattice.}

{\bf Proof.} Let $\{e_i\}$ be a set of idempotents of $A$.
Obviously $e_i\leq 1$ for every index $i$. If there does not exist
a minimal idempotent $q_1$ such that $e_i\leq 1-q_1$ for every
index $i$ then $\sup_i e_i=1$, otherwise, if there does not exist
a minimal idempotent $q_2$ such that $e_i\leq 1-q_1-q_2$ for every
index $i$ then $\sup_i e_i=1-q_1$ and so on. Here $q_1$,
$q_2$,$\dots$ are pairwise orthogonal idempotents. Hence the
sequence $q_1$, $q_2$,$\dots$ is finite and $\{e_i\}$ has a least
upper bound in $A$. This gives us
$$
\inf_i e_i=1-\sup_i(1 - e_i).
$$
Thus set of all idempotents of $A$ is a complete lattice.
$\triangleright$

{\bf Corollary 3.6.} {\it Every finite dimensional unital Jordan
algebra without nilpotent elements is a Baer Jordan algebra.}

By corollary 3.6 and lemma 2.4 the following theorem is valid.

{\bf Theorem 3.7.} {\it Let $A$ be a finite dimensional Jordan
algebra $A$ with identity. Then the following statements are
equivalent:

1) $A$ has no nilpotent elements,

2) $A$ is a Baer Jordan algebra.}

{\bf Theorem 3.8.} {\it Let $A$ be a finite dimensional unital
Jordan algebra $A$ without nilpotent elements. Suppose that $A$
contains $n\geq 3$ pairwise orthogonal and strongly connected
minimal idempotents with sum 1. Then $A$ is isomorphic to one of
$H_n({\Bbb R})$, $H_n({\Bbb C})$, $H_n({\Bbb H})$ or, if $n=3$,
$H_3({\Bbb O})$.}

{\bf Proof.} The coordinatization theorem (2.8.9 in [4])
tells us that $A\cong H_n(R)$, for some real $*$ algebra $R$.

Let $a\in R$ be an arbitrary element. Then $ae_{ij}+a^*e_{ji}\in
\{e_iAe_j\}$. Hence
$$
(ae_{ij}+a^*e_{ji})^2=aa^*e_{ii}+a^*ae_{jj}=\lambda(e_{ii}+e_{jj})
$$
and $aa^*e_{ii}=\lambda e_{ii}$ by (i) of lemma 3.1. Therefore
$aa^*=\lambda 1$, $\lambda\in {\Bbb R}$ and $R_{sa}={\Bbb R}1$.

Suppose that $a\in R$ is not zero. Let $x=ae_{12}+a^*e_{21}$. Then
$$
0\neq x^2=(aa^*)e_{11}+(a^*a)e_{22}.
$$
Since the Jordan algebra $A$ is without nilpotent elements we have
either $aa^*$ or $a^*a$ is nonzero.

However, since $a+a^*\in R_{sa}={\Bbb R}1$, $a$ and $a^*$ commute,
so $a^*a\neq 0$. By 2.7.8 in [4] $R$ is $*$ isomorphic to
either $\Bbb R$, $\Bbb C$, $\Bbb H$ or, if $n=3$, to $\Bbb O$.
$\triangleright$

Since every finite dimensional Baer Jordan algebra has no nonzero
nilpotent elements we have the following corollary.

{\bf Corollary 3.9.} {\it Every finite dimensional Baer Jordan
algebra $A$ with $n\geq 3$ pairwise orthogonal and strongly
connected minimal idempotents with sum 1 is isomorphic to one of
$H_n({\Bbb R})$, $H_n({\Bbb C})$, $H_n({\Bbb H})$ or, if $n=3$,
$H_3({\Bbb O})$.}

{\bf Theorem 3.10.} {\it Let $A$ be a finite dimensional (with
respect to field ${\Bbb R}$), formally real, unital Jordan algebra
which also contains two connected minimal idempotents with sum
$1$. Then

1) if $\dim(A)=3$, then $A\cong H_2({\Bbb R})$,

2) if $\dim(A)=4$, then $A\cong {\Bbb R}\oplus H$, where $H$ is a
two dimensional Hilbert space.}

{\it Proof.} Let $p$ and $q$ be minimal idempotents in such
algebra $A$, with $p+q=1$. Clearly, $p$ and $q$ are orthogonal.
The following Peirce decomposition holds
$$
A={\Bbb R}p\oplus {\Bbb R}q \oplus A_{12}.
$$
where $A_{12}=\{pAq\}$.

1): It is clear that, if $A$ is three dimensional then there
exists an element $s\in A_{12}$ such that
$$
A_{12}={\Bbb R}s, s^2=p+q.
$$
Then $U_sp=2s(sp)-s^2p=s^2-p=p+q-p=q$. Similarly $U_sq=p$. Hence
$A\cong H_2({\Bbb R})$.

2): Now suppose that $A$ is four dimensional. In this case the
vector space $A_{12}$ is two dimensional. Let $\{e_1, e_2\}$ be a
basis of $A_{12}$. By (v) of lemma 2.9.4 in [4]
$e_i^2=\lambda_i (p+q)$, $i=1,2$ for some positive numbers
$\lambda_1$ and $\lambda_2$ in ${\Bbb R}$. Without loss of
generality we assume that $e_i^2=p+q$, $i=1,2$. Let
$$
s_1=p-q, s_2=e_1+e_2, s_3=e_1-e_2.
$$
Then $s_is_j=0$ if $i\neq j$. Also we may assume $s_i^2=1$,
$i=1,2,3$. Indeed, we have
$$
e_1+e_2\in A_{12}, (e_1+e_2)^2=\lambda (p+q),
$$
for some $\lambda\in {\Bbb R}_+$. At the same time,
$$
(e_1+e_2)^2=e_1^2+2e_1e_2+e_2^2=2(p+q)+2(\alpha p+\beta q),
$$
for some real numbers $\alpha$, $\beta\in {\Bbb R}$ such that
$e_1e_2=\alpha p+\beta q$. Hence
$\alpha=\beta=\frac{\lambda-2}{2}$, $\lambda=2\alpha+2$. Since
$\lambda>0$ we have $\alpha>-1$. Hence we may assume that
$$
s_2=\frac{1}{\sqrt{2\alpha+2}}(e_1+e_2).
$$
Similarly we have
$$
(e_1-e_2)^2=(2-2\alpha)(p+q), 1>\alpha.
$$
Hence $-1<\alpha<1$ and we may assume that
$$
s_3=\frac{1}{\sqrt{2-2\alpha}}(e_1-e_2).
$$
Therefore $\{s_1, s_2, s_3\}$ is a spin system and $A$ is a spin
factor. $\triangleright$

\end{document}